\documentclass[11pt]{amsart}

\usepackage{amstext}
\usepackage{amsthm}
\usepackage{amsopn}

\usepackage{amsfonts}
\usepackage{amsmath}
\usepackage{latexsym}
\usepackage{amscd}
\usepackage{amssymb}
\usepackage{amsmath}

\swapnumbers
\newtheorem{thm}[equation]{Theorem}
\newtheorem{prop}[equation]{Proposition}
\newtheorem{lem}[equation]{Lemma}
\newtheorem{cor}[equation]{Corollary}

\theoremstyle{definition}

\newtheorem{remark}[equation]{Remark}

\numberwithin{equation}{section}

\def\limind{\mathop{\oalign{lim\cr
\hidewidth$\longrightarrow$\hidewidth\cr}}}

\newcommand{\bbA}{{\mathbb A}}
\newcommand{\bbG}{{\mathbb G}}

\newcommand{\bbZ}{{\mathbb Z}}

\newcommand{\bbQ}{{\mathbb Q}}

\newcommand{\Aut}{\operatorname{Aut}}
\newcommand{\Spec}{\operatorname{Spec}}

\newcommand{\rad}{\operatorname{rad}}

\newcommand{\Pic}{{\operatorname{Pic}}}
\newcommand{\simlgr}{\buildrel \sim \over \lra}

\newcommand{\rank}{\operatorname{rank}}

\newcommand{\Char}{\operatorname{char}} %% \char is already a command

\newcommand{\Gal}{\operatorname{Gal}}

\newcommand{\lra}{\longrightarrow}

\newcommand{\cal}{\mathcal}

\begin{document}
\title[$G$-torsors]{Reduction of structure for torsors
over semilocal rings}

\author{V. Chernousov}
\address{Department of Mathematics, University of Alberta,
    Edmonton, Alberta T6G 2G1, Canada}
\thanks{ V. Chernousov was partially supported by the Canada 
Research Chairs Program and an NSERC research grant} 
\email{chernous@math.ualberta.ca}

\author{P. Gille}
\address{UMR 8552 du CNRS, DMA, Ecole normale
sup\'erieure, F-75005 Paris, France}
\email{Philippe.Gille@ens.fr}

\author{Z. Reichstein}
\address{Department of Mathematics, University of British Columbia,
       Vancouver, BC V6T 1Z2, Canada}
\email{reichst\@@math.ubc.ca}
\urladdr{www.math.ubc.ca/$\stackrel{\sim}{\phantom{.}}$reichst}
\thanks{Z. Reichstein was partially supported by an NSERC research grant}


%%%%%%%%%%%%%%%%%%%%%%%%%%%%%%%%%%%
% 14L30 Group actions on varieties and schemes
% 14E Birational geometry,
% 14E20 Coverings
% 11E72 Galois cohomology of linear algebraic groups
%%%%%%%%%%%%%%%%%%%%%%%%%%%%%%%%%%%
% Better numbers for this paper
%
% 14L15 group schemes
% 20G10 Cohomology theory of algebraic groups
%
% 11E72 is  still OK
%%%%%%%%%%%%%%%%%%%%%%%%%%%%%%%%%%
%%%%%%%%%%%%%%%%%%%%%%%%%%%%%%%%%%%


\begin{abstract}
Let $G$ be a reductive affine group scheme defined over
a semilocal ring $k$.
Assume that either $G$ is semisimple or $k$ is normal and
noetherian.  We show that $G$ has a finite $k$-subgroup
$S$ such that the natural map $H^1(R, S) \to H^1(R, G)$
is surjective for every semilocal ring $R$ containing $k$.
In other words, $G$-torsors over $\Spec(R)$ admit reduction
of structure to $S$.  We also show that the
natural map $H^1(X, S) \to H^1(X, G)$ is surjective
in several other contexts, under suitable assumptions on the
base ring $k$, the scheme $X/k$ and the group scheme
$G/k$.  These results have already been used
to study loop algebras and essential dimension
of connected algebraic groups in prime characteristic.
Additional applications are presented at the end of this paper.
% In this paper we give several new applications.
\end{abstract}

\maketitle
 
\tableofcontents

\section{Introduction}

Let $G$ be a linear algebraic group defined over a field $k$.
In~\cite{cgr} we showed that, under mild
assumptions on $G$ and $k$, $G$ has a finite $k$-subgroup $S$ such that
every $G$-torsor over a field $K/k$ admits reduction of structure
to $S$, i.e.,~the natural map of Galois cohomology sets
$H^1(K, S) \to H^1(K, G)$ is surjective. In several subsequent
applications, a more general version of this result was
needed, with the field $k$ replaced by a base ring, the group
$G$ by a reductive group scheme over $k$ and the field
$K/k$ by a $k$-scheme $X$. The goal of this paper
is to extend the main result of \cite{cgr} to this more
general setting.  

{\em All schemes in this paper will be assumed to be locally noetherian.}
Of particular interest to us will
be $k$-schemes $X$ satisfying the following condition:
\begin{equation} \label{e.C}
\text{$\Pic(X') = 0$ for every generalized Galois cover $X'/X$.}
\end{equation}
Here by a generalized Galois cover $X' \to X$ we mean a $\Gamma$-torsor,
for some twisted finite constant group scheme $\Gamma$ defined 
over $X$. In other words, $\Gamma = {_aC}$, where $C$ is 
a finite constant group scheme over $X$ and $[a] \in H^1(X, \Aut(C))$. 
(The term ``Galois cover" is usually reserved for
the case where $\Gamma$ is itself a finite constant group scheme.) 
The class of schemes satisfying
condition~\eqref{e.C} includes, in particular, affine schemes
of the form $X = \Spec(R)$, where $R$ is a semilocal
ring containing $k$.  If $K$ is a $k$-field of
characteristic 0, we can also take $R$ to be
a polynomial ring $K[x_1, \dots, x_n]$
(see~\S\ref{sect.affine}) or a Laurent polynomial ring
$K[x_1^{\pm 1}, \dots, x_n^{\pm 1}]$ (see Remark~\ref{rem.loop}).

We are now ready to state the main results of this paper. Recall
that an $X$-group $T$ of multiplicative type is called {\em
isotrivial} if $T \times_X X'$ is split for some finite \'etale 
surjective map $X' \to X$. For the definition and basic properties
of groups of multiplicative type, we refer the reader
to~\cite[X]{SGA3}.

\begin{thm} \label{thm1}
Let $k$ be a commutative base ring and $G$ be a smooth affine
group scheme over $k$ whose connected component $G^0$ is reductive.
Assume further that one of the following holds:

\medskip

(a) $k$ is an algebraically closed field, or

\medskip

(b) $k=\bbZ$, $G^0$ is a split Chevalley group, and
the order of the Weyl group of the geometric fiber 
$G_{\overline{s}}$ is independent of $s \in \Spec(\bbZ)$, or 

\medskip

(c) $k$ is a semilocal ring, $G$ is connected, and
the radical torus $\rad(G)$ is isotrivial.

\medskip

\noindent
Then there exist a maximal torus $T \subset G$ defined over $k$ and 
a finite $k$-subgroup $S \subset N_G(T)$,
such that

\smallskip
\begin{enumerate}
\item $S$ is an extension of a twisted constant group
scheme by a finite $k$-group of multiplicative type,

\smallskip
\item the natural map $ H^1(X, S) \lra H^1(X, N_G(T)) $
is surjective for any scheme $X/k$ satisfying condition~\eqref{e.C}.
\end{enumerate}
\end{thm}

Of course, if $G$ is connected then (a) is a special case of (c).
Note also that in case $(b)$ we can take $G$ to be
the automorphism group $\Aut(G_0)$ of some
semisimple Chevalley group scheme $G_0$. In this case
the cohomology set $H^1(X, G)$ classifies the
semisimple group schemes over $X$,
which are \'etale locally isomorphic to $G_0 \times_{\bbZ} X$.

\smallskip

Combining Theorem~\ref{thm1}(c) with Grothendieck's existence
theorem for maximal tori (reproduced as Theorem~\ref{thm3}
below), we obtain the following stronger result in case (c);
cf. \S\ref{sect.thm2}.

\begin{thm} \label{thm2}
Let $k$, $G$ and $S$ be as in Theorem~\ref{thm1}(c).
Then the map $H^1(R, S) \to H^1(R, G)$ is surjective 
for any semilocal ring $R/k$.
\end{thm}

Note that the assumption on the radical of $G$ is superfluous
if $G$ is a semisimple group scheme or if $k$ is normal and
noetherian, because all tori defined over such rings
are isotrivial; see~\cite[X.5.16]{SGA3}.

The symbol $H^1(X, G)$ in the statements of Theorems~\ref{thm1}
and~\ref{thm2} denotes the flat cohomology set, which
classifies $G$-torsors over $X$; see~\S\ref{sect.prel}.
If $G$ is smooth then every $G$-torsor over $X$ is also smooth
and is trivialized by an \'etale covering \cite[III.4]{M}. 
So in this case the natural map $H^1_{\acute et}(X, G)
\to H^1(X, G)$ is bijective, and we may replace $H^1(X, G)$ by
$H^1_{\acute et}(X, G)$.

In particular, suppose that $k$ is an algebraically
closed field and $G/k$ and $S/k$ are as in Theorem \ref{thm1}(a).
If $K$ is a perfect field containing $k$ and $K_s$
is the separable closure of $K$ then
\[ H^1(K,G) = H^1(\Gal(K_s/K),G(K_s)) \]
and 
\[ H^1(K,S) = H^1(\Gal(K_s/K),S(K_s)) \, . \]
In other words, in this situation the flat cohomology
sets appearing in the statement
of Theorem~\ref{thm1}(a) can be replaced by Galois cohomology.
Moreover, since $S$ is finite, $S(k)=S(K_s)$ (with
$\Gal(K_s/K)$ acting trivially on both sides) and hence,
$$
H^1(\Gal(K_s/K),S(k)) =  H^1(\Gal(K_s/K),S(K_s)) \, .
$$
Thus in this setting Theorem~\ref{thm1}(a) implies the following
characteristic-free result about Galois cohomology.
The assertion about $|S| := \dim_k \, k[S]$ is immediate 
from the construction of $S$ in \S\ref{sect.proof-of-1ab}.

\begin{cor} \label{cor1}
Let $G$ be a linear algebraic group defined over an
algebraically closed field $k$, whose connected component
$G^0$ is reductive.  Then there exists a finite
$k$-subgroup $S$ of $G$, such that
every prime factor of $|S|$ divides the order of the
Weyl group $W(G)$, and the map
$$
H^1(\Gal(K_s/K),S(k)) \to H^1(\Gal(K_s/K),G(K_s))
$$
is surjective for any perfect field $K/k$.
\qed
\end{cor}

Corollary~\ref{cor1} generalizes~\cite[Theorem 1.1(a)]{cgr}, which
yields the same conclusion if $\Char(k) = 0$. This corollary has been
used to study essential dimension of connected algebraic groups in
positive characteristic in~\cite{GR}. An application of
Theorem~\ref{thm1} to the study of loop algebras can be found
in~\cite{GP2}; cf. Remark~\ref{rem.loop}. We will give additional 
applications in~\S \S\,\ref{sect.thm2} and~\ref{sect.affine}. 
A further application of Theorem~\ref{thm1}
will appear in the forthcoming paper~\cite{cgp}.

We are grateful to  M. Florence for pointing out a mistake
in the proof of \cite[Theorem~1.1(a)]{cgr}. This mistake is
corrected in the course of the proof of~Lemma~\ref{lem2}.
For details, see Remark~\ref{rem.correction}.

\section{Preliminaries}
\label{sect.prel}

We begin by recalling some known facts about affine group
schemes $G$ of finite type over an arbitrary  base scheme $X$.

{\em A pseudo $G$-torsor}
({\em formellement principal homog\`ene}, in~\cite{SGA3})
$E$ over $X$ is an $X$-scheme
equipped with a right action of $G$ such that the mapping $E
\times_X G \to E \times_X E$ given by $(x,g) \mapsto (x, x.g)$
is an isomorphism; 
see~\cite[{\uppercase\expandafter{\romannumeral 4}}.5.1]{SGA3}. 
A pseudo $G$-torsor $E$ is {\em a $G$-torsor}
({\em fibr\'e principal homog\`ene})
if it is locally trivial in the
fppf topology, i.e., if
there exists a faithfully flat morphism $X' \to X$, locally
of finite type, such that $E \times_X X' \cong G \times_X X'$.
Here, as usual, the acronym fppf stands for
``fid\`element plate de pr\'esentation finie" or
``faithfully flat and finitely presented".

For such a covering  $X' \to X$, we define
\[ Z^1(X'/X , G) := 
 \bigl\{  g \in G(X'  \times_X X') \, \mid
\, p_{1,2}^*(g) p_{2,3}^*(g)=p_{1,3}^*(g) \bigr\} \] 
% \in G(X' \times_X X' \times_X X') \bigr\}  \]
and
$$
H^1(X'/X,G):= Z^1(X'/X,G) / G(X'),
$$
where  $G(X')$ acts on  $Z^1(X'/X,G)$ by $g \cdot z= p_1^*(g) \, z \,
p_2^*(g)^{-1}$;
see~\cite[Chapter {\uppercase\expandafter{\romannumeral 3}}]{K}.
Here $p_{i,j} \colon X' \times_X X' \times_X X' \to X' \times_X X'$
is the projection \[ (x_1', x_2', x_3') \to (x_i', x_j') \, , \]
and 
$p_{1,2}^*(g)$, $p_{2,3}^*(g)$, $p_{1,3}^*(g)$ are viewed as elements of
% the identity
% $p_{1,2}^*(g) p_{2,3}^*(g)=p_{1,3}^*(g)$ is required to 
% hold in $G(X' \times_X X' \times_X X')$.
$G(X' \times_X X' \times_X X')$.
The pointed set $H^1(X'/X,G)$ classifies  $G$-torsors over
$X$ which are trivialized by the base change $X'/X$,
i.e.,~$G$-torsors $E$ satisfying 
\[ E \times_X X'
{\buildrel\sim\over{\longrightarrow}} G \times_X X' \, ; \]
see \cite[{\uppercase\expandafter{\romannumeral 3}}.4, page 120]{M}.  
We now define
 $$ H^1(X,G):=\limind_{X'} H^1(X'/X,G) \, , $$
where the limit
is taken over all coverings $X'/X$ in the $fppf$ topology.
The pointed set $H^1(X,G)$ classifies $G$-torsors over $X$.

If $P$ is a $G$-torsor over $X$, we denote by ${^P}G$ the associated
twisted $X$-group scheme; it is the twisted inner form of $G$ and can
be defined as the scheme of $G$-automorphisms of $P$.
We then have a canonical bijection (the ``torsion'' map)
$$
% H^1(X,G) \longrightarrow H^1(X,{^P}G)
H^1(X,G) \xrightarrow{\simeq} H^1(X,{^P}G)
$$
mapping  a $G$-torsor $Q$ to the  scheme
$\underline{\mathrm{Isom}}_G(P,Q)$ of  $G$-isomorphisms of  
$P$ into $Q$; see~\cite[III.2.6]{Gir}. In particular, 
the torsion map takes $P$ to the trivial ${^P}G$-torsor.

We say that $G$ is {\em connected} if the fiber $G_x$ is connected for
any point $x \in X$. Here we view $G_x = G \times_X \Spec(\kappa(x))$ 
as an algebraic group over the residue field $\kappa(x)$ of $x$.
If $G/X$ is smooth, then $G$ contains a unique maximal open connected
normal subgroup defined over $X$; \cite[VI$_B$, Thm. 3.10]{SGA3}.
As usual, we will denote this subgroup by $G^0/X$ and refer 
to it as {\em the connected component} of $G$.  Note that $G^0$ 
is smooth over $X$ and it is a closed subgroup of $G$; 
in particular, it is affine over $X$.

We say that $G/X$ is {\em reductive} if it is smooth and all
of its geometric fibers $G_{\overline{x}}$ are (connected) reductive
groups~\cite[XIX.2.7]{SGA3}.  A subgroup $T/X$ of $G/X$ is a
{\em maximal torus} if it is an $X$-torus and all of its geometric 
fibers are maximal tori \cite[XII.1.3]{SGA3}. The {\em radical torus}
$\rad(G)$ of $G$ is the unique maximal torus of the center
of $G$~\cite[XXII.4.3.6]{SGA3}.

Similarly a subgroup $B/X$ of a reductive
group scheme $G/X$ is a {\em Borel subgroup} if it is smooth
and finitely presented and all of its geometric fibers are Borel
subgroups \cite[XXII.5.2.3]{SGA3}.

We refer to \cite[XXII.1]{SGA3} for the definitions of split
group schemes and to~\cite[XXIV.3]{SGA3} for the definition
of the Dynkin scheme of $G$ and quasi-split reductive
group schemes.

Let $G$ be a split adjoint semisimple group over $X$,
$T$ a maximal split torus in $G$ defined over $X$,
$B$ a Borel subgroup containing $T$ and $D/X$
the corresponding Dynkin scheme of $G$.
Following~\cite[XXIV.3.5]{SGA3}, we will denote
the group scheme representing the functor of
automorphisms of $D$ (as a Dynkin scheme) by
$\Aut_{Dyn}(D)$.
By \cite[XXIV, Th\'eor\`eme 1.3 and 3.6]{SGA3}
$$
\Aut(G)=G\rtimes \Aut_{Dyn}(D) \, .
$$
Moreover, there exists a canonical splitting
$$ h \colon \Aut_{Dyn}(D) \to \Aut(G)$$
such that the image of $h$ preserves
$T$ and $B$. Every quasi-split adjoint group scheme
$G'$ of the same type as $G$ is
$X$-isomorphic to the twist ${_{h_*(a)}G}$ of $G$ for
some cocycle $$ a\in Z^1_{\acute et}(X,\Aut_{Dyn} (D)) \, . $$

\section{A first step towards the proof of Theorem~\ref{thm1}}
\label{sect.proof-of-0}

The purpose of this section is to prove the following proposition,
which will play a key role in the proof of Theorem~\ref{thm1}.

\begin{prop} \label{prop0}
Let $k$ be a commutative base ring and
\[ 1 \to T \to N \xrightarrow{p} W \to 1 \]
be an exact sequence of smooth group schemes defined over $k$,
where $T$ is an isotrivial torus, split by a Galois extension $k'/k$
of degree $d$, and $W$ is a twisted finite constant group
of order $n$.  Suppose $N$ has a finite $k$-subgroup $S'$ such that
$p(S') = W$.  Then there exists a finite $k$-subgroup $S \subset G$
containing $S'$ such that the natural map $H^1(X, S) \to H^1(X, N)$ 
is surjective for any $k$-scheme $X$ satisfying condition~\eqref{e.C}.

Moreover, we can take $S$ to be the subgroup of $N$ generated by
$S'$ and $\phi_{m}^{-1}(S' \cap T)$, where $m = nd$ and
$\phi_{m} \colon T \to T$ is the map taking $t \in T$ to $t^{m}$.
\end{prop}

\begin{proof}
Denote by $q \colon S \to W$ and $q' \colon S' \to W$ the restrictions
of the projection $p \colon N \to W$ to $S$ and $S'$, and 
by $\mu = S \cap T$ and $\mu' = S' \cap T$ the kernels of these maps,
respectively.  Let $X$ be a $k$-scheme satisfying condition~\eqref{e.C} 
on Picard groups.  We will prove the surjectivity 
of $H^1(X, S) \rightarrow H^1(X,N)$ fiberwise, 
with respect to the mapping  $p_*: H^1(X,N) \rightarrow H^1(X,W)$ 
induced by $p$. Fix $[b] \in
H^1(X,N)$; our goal is to show that $[b]$ lifts to
$H^1(X, S)$. 

\begin{lem} \label{lem2}
Let $[a]=p_*([b]) \in H^1(X,W)$. Then $$[a]\in {\rm Im}\,(
H^1(X,S) \buildrel q_* \over \to H^1(X,W)).$$
\end{lem}

\begin{proof}[Proof of Lemma~\ref{lem2}]
The obstruction to lifting $[a]$ to $H^1(X, S)$ is the class
$$
\Delta([a]) \in  H^2(X, {_a\mu}),
$$
where ${_a\mu}$ denotes the group $\mu$ twisted by the torsor $a$
\cite[IV.4.2.8]{Gir}. We now use the commutative diagram 
\begin{equation} \label{diagr.lem2}
\begin{CD}
1 @>>> \mu' @>>> S' @>{q'}>> W @>>> 1 \\
&&  \cap &&  \cap && \mid \, \mid \\
1 @>>> \mu @>>> S @>{q}>> W @>>> 1 \\
&&  \cap &&  \cap && \mid \, \mid \\
1 @>>> T @>>> N @>{p}>> W @>>> 1 \\
\end{CD}
\end{equation}
with exact rows and the functoriality of the obstruction 
$\Delta([a])$. If $\Delta'([a]) \in H^2(X, \, _a\mu')$ 
is the obstruction to lifting $[a]$ to
$H^1(X, S')$, via $q'_*: H^1(X, S') \rightarrow
H^1(X,W)$, then $\Delta([a])$ is the image of $\Delta'([a])$
under the natural map $H^2(X, \, {_a\mu'}) \to H^2(X, \,
{_a\mu})$.

The commutative diagram 
\begin{equation} \label{e.S}
\begin{CD}
1 @>>> {_a\mu'} @>>> {_aT} @>{t \to \;  (t \; \text{mod} \, \mu') }>> {_a(T/\mu')} @>>> 1 \\
&&  \cap && \mid \, \mid   && @VVV \\
1 @>>> {_a\mu} @>>> {_aT} @>{t \to \; (t^m \; \text{mod} \, \mu')}>> {_a(T/\mu')} @>>> 1 \\
\end{CD}
\end{equation}
with exact rows gives rise to the commutative exact diagram
$$
\begin{CD}
H^1(X, \, {_a(T/\mu')}) @>>> H^2(X, \, _a\mu') @>>> H^2(X, \, {_aT})   \\
 @V{\times m}VV  @VVV      \mid \, \mid    \\
H^1(X, \, {_a(T/\mu')}) @>>> H^2(X, \, {_a\mu}) @>>> H^2(X, \, {_aT}) 
\end{CD}
$$
which we will now analyse. Recall that the middle vertical map
sends $\Delta'([a]) \in H^2(X, \, _a\mu')$ 
to $\Delta([a]) \in  H^2(X, \, {_a\mu})$. Since we 
are given that $[a]$ lifts to $[b] \in H^1(X, N)$, 
we have
$$
\Delta'([a]) \in  \ker\Bigl( H^2(X, {_a\mu'})\to
H^2(X, {_aT}) \Bigr)  $$
and thus
$$\Delta'([a])
\in  {\rm Im}\,\Bigl( H^1(X, {_a(T/\mu')})\to H^2(X,
{_a\mu'}) \Bigr). $$ 
In order to prove the lemma (i.e., to
prove that $\Delta([a]) = 0$), it now
suffices to show that the vertical map
\begin{equation} \label{e.times-m}
\begin{CD}
H^1(X, \, {_a(T/\mu')})   \\
 @V{\times m}VV    \\
H^1(X, \, {_a(T/\mu')})
\end{CD}
\end{equation}
in the above diagram is trivial.

If $p:X'\to X$ is a cover (i.e., a finite \'etale map) of degree $m$
and $H$ is a commutative affine $X$-group scheme, we will denote 
the trace morphism by $N_{X'/X}:{\rm R}_{X'/X}(H)\to H$; 
cf.~\cite[0.4]{CTS2}.  If $p$ has degree $m$, the composition
$$
H \longrightarrow {\rm R}_{X'/X}(H)
\xrightarrow{N_{X'/X}} H
$$
of $N_{X'/X}$ with the natural map
$H \to {\rm R}_{X'/X}(H)$ is multiplication by $m$.
%%%%%%%%%%%%%%%%%%%%%%%%%%%%%%%%%%%%%%%%%%%%

Now let $Y \to X$ be the $W$-torsor associated to $a$ and apply
the above facts to the generalized Galois covering  
$X'=Y \times_k k' \to X$
of degree $m = nd$, with $H={_a(T/\mu')}$. Note that
this covering trivializes $a$ and splits $T$. 
The map~\eqref{e.times-m} can be decomposed as
$$
H^1(X,{_aT})\longrightarrow H^1(X,{\rm R}_{X'/X}({_a(T/\mu')}))
\longrightarrow H^1(X,{_a(T/\mu')}).
$$
Shapiro's lemma and condition~\eqref{e.C} imply that
$$
H^1(X,{\rm R}_{X'/X}({_a(T/\mu')}))=H^1(X',{_a(T/\mu')})=
\Pic\,(X')^{\rank(T)}=0 \, .
$$
Hence the map~\eqref{e.times-m} is trivial, as claimed. The proof
of Lemma~\ref{lem2} is now complete.
\end{proof}

We are now ready to finish the proof of Proposition~\ref{prop0}. Let
$[c] \in H^1(X,S)$ be such that $q_*([c])=[a]$. The bottom two rows
of~\eqref{diagr.lem2} give rise to the diagram
\[ \begin{CD}
H^1(X, {_c \mu}) @>>> q_*^{-1}(a) \subset H^1(X, S) \\
 f @VVV  @VVV  \\
H^1(X, {_c T}) @>>> p_*^{-1}(a) \subset H^1(X, N)
\end{CD} \]
where the horizontal arrows are the ``torsion" maps (see \S 2). Recall
that our goal is to show that $[b] \in p_*^{-1}([a]) \subset
H^1(X, N)$ lies in the image of $H^1(X, S)$. If $X = \Spec(K)$
for some field $K/k$ then a twisting argument~\cite[I.5.5]{serre-gc}
shows that the map $$ H^1(K, {_{c}T}) \rightarrow p_*^{-1}([a])$$ is
surjective. The same twisting argument goes through
for any $k$-scheme $X$~\cite[III.3.2.4]{Gir}; in this case
we can also conclude that
the map $$ H^1(X, {_{c}T}) \rightarrow p_*^{-1}([a])$$ is
surjective.  Thus it suffices to prove that the vertical map
$f$ in the above diagram
is surjective as well.  The exact sequence
\begin{equation} \label{e.cT}
\begin{CD} 
 1 @>>> _c \mu @>>> _c T @>{t \to \; (t^m \; \text{mod} \, \mu')}>> 
_c {T/\mu'} @>>> 1 
\end{CD} 
\end{equation}
gives rise to the exact sequence 
\[ \begin{CD}
H^1(X, {_c \mu}) @>f>> H^1(X, \, _cT) @>>> H^1(X, \, _c(T/\mu)) \, .  
\end{CD} \]
It thus remains to show that the map
\begin{equation} \label{e.prop0-final}
H^1(X,{_c T}) \to H^1(X,{_{c}(T/\mu)})
\end{equation}
in this sequence is trivial. Indeed, since the group homomorphism
\[ _c T \to \, _c(T/\mu') \] in~\eqref{e.cT} factors through 
\[ \times m \colon _c(T/\mu') \xrightarrow{ \; \; \;  \times m \; \; \;}  
\, _c(T/\mu') \, , \]
the map \eqref{e.prop0-final} factors through 
\[ \times m : H^1(X, \, _c(T/\mu')) \xrightarrow{\quad}  H^1(X, \, _c(T/\mu')) \]
which we showed to be trivial
at the end of the proof of Lemma~\ref{lem2}.
We conclude that the map~\eqref{e.prop0-final} is trivial, as claimed.
This completes the proof of Proposition~\ref{prop0}.
\end{proof}

\begin{remark} \label{rem.prop0}
Let $k$ be a ring, $T$ is a maximal $k$-torus in an affine algebraic 
$k$-group $G$ and $N = N_G(T)$. This is a natural setting, 
where Proposition~\ref{prop0} can be applied. However, 
it is not a priori clear for which $G$ one can construct
a finite group $S'$ as in Proposition~\ref{prop0}.
In fact, it is not even clear in general which affine $k$-groups 
$G$ contain a maximal $k$-torus $T$. If we can find 
a maximal $k$-torus $T \subset G$ and a finite $k$-subgroup
$S' \subset N = N_G(T)$ with desired properties, we would also
like to know under what circumstances one can conclude that the map
$H^1(X, S) \to H^1(X, G)$ is surjective. In the sequel we will 
give partial answers to these questions, under additional 
assumptions on $k$.
\end{remark}

\section{Proof of Theorem \ref{thm1}(a) and (b)}
\label{sect.proof-of-1ab}

(a) Let $T$ be a maximal $k$-torus of $G$ and $N = N_G(T)$.  Since
\[ N^0 \subset \bigl( N_{G^0}(T) \bigr)^0 =T \subset N \, , \]
we have $N^0=T$. Hence, $N$ is smooth and $W$ is a finite constant 
group. Let $p \colon N \to W = N/T$ be the natural projection. 
By Proposition~\ref{prop0} 
it suffices to construct a finite $k$-subgroup $S' \subset N$
such that $p(S') = W$. In fact, we will construct $S'$ so that
$\mu'$ be the $n$-torsion subgroup of $T$, where $n = |W|$. 

Consider the exact sequences
\[
1\to T \to N \stackrel{p}{\to} W \to 1 \quad \quad  \text{and} \quad
\quad 1 \to  \mu' \to T \buildrel \times n \over \to T \to 1 \, . \]
According to \cite[II.2, Proposition 2.3]{DG} (cf. also
\cite[XVII, App. I.3.1, page 622]{SGA3}),
extensions of $W$ by $T$ are classified by the Hochschild cohomology
group $H^2_0(W,T)$. Since $W$ is a constant group scheme,
$H^2_0(W,T)$ is isomorphic to the usual cohomology group
$H^2(W,T(k))$; see \cite[III.6.4, Proposition 4.2]{DG}. Thus 
the first sequence yields a class in $H^2(W,T(k))$. Since $n \cdot
H^2(W,T(k))=0$, the second sequence tells us that this class comes
from $H^2(W, \mu'(k))$. In other words, there is
an extension $S' \subset N$ of $W$ by $\mu'$ such 
that $N$ is the push-out of $S'$ by the morphism $\mu' \hookrightarrow T$.
This completes the construction of $S'$.

\smallskip
(b) Let $T$ be a maximal split torus of $G$ defined over $\bbZ$. 
Note that $W = N(T)/T$ is a constant finite group 
scheme; this follows from the fact that $W$ is representable by
a $\bbZ$-group scheme which is finite \'etale \cite[XII.2.1.b]{SGA3}.

It remains to construct a finite subgroup $S' \subset N$
which surjects onto $W$; the desired conclusion will then follow
from Proposition~\ref{prop0}.

Our construction of $S'$ will be based on schematic 
adherence, which associates to a closed $\bbQ$--subscheme
$V \subset G_{\bbQ}$ its Zariski closure $\widetilde V$ in $G_{\bbZ}$.
Schematic adherence induces a one-to-one correspondence between
$\bbQ$--subchemes of $G_{\bbQ}$ and flat closed $\bbZ$--subchemes of
$G_{\bbZ}$~\cite[I.2.6]{BT}.  In particular, it maps $\bbQ$--subgroups 
of $G_{\bbQ}$ into flat $\bbZ$--group subschemes 
of $G$~\cite[I.2.7]{BT} (see also \cite[\S 3]{GM}).

Let $n=|W|$ and $T=\bbG_m^r$, where $r$ is the rank
of $G$.  As pointed out by Tits \cite{tits2}, the fact 
that $H^1(\bbZ,T)=\Pic(\bbZ)^r=0$ implies that the sequence
$$
0 \to T(\bbZ) \to N(\bbZ) \to W \to 1.
$$
is exact.  Since $T(\bbZ)= \{ \pm 1 \}^r$, 
$N(\bbZ)$ is a finite group. View $N(\bbZ)$
as a finite constant $\bbQ$--subgroup of $G_{\bbQ}$
and let $S'$ be its schematic adherence 
in $N/ {\bbZ}$. Then $S'$ is a finite
flat $\bbZ$--subgroup scheme of $N$.
Since $N(\bbZ)$ surjects onto $W$, so does $S'$.
\qed

%  contains $B(\
%  equipped with a natural
% $\bbZ$--morphism
% $\widetilde p:  \to W$.
% Since ${\widetilde  N}(\bbZ)= G(\bbZ) \cap N(\bbQ)= N(\bbZ)$,
% $\widetilde p$ is a surjective morphism of group schemes.
% We claim that $\ker(\widetilde p)= {_2T}$.
%% I do not know to prove directly that $\ker(\widetilde p)$ is flat.
% That kernel is a subscheme of the schematic adherence
% of ${\widetilde  N}_{\bbQ} \cap T_{\bbQ}$ in $T$ which is
% the schematic adherence of the constant group
% of $T(\bbZ) = N(\bbZ) \cap T_{\bbQ}$  in $T$, namely
% ${_2T}$. Hence $\ker(\widetilde p) \subset {_2T}$ and 
% the converse inclusion is obvious.  Thus we have the
% following compatible exact sequences of flat group
% $\bbZ$-schemes
% $$
% \begin{CD}
% 0 @>>> \mu_2^r @>>> \widetilde N @>{\widetilde p}>> W @>>> 1 \\
% &&@VVV @VVV  \mid \mid \\
% 0 @>>> T  @>>> N @>>> W @>>> 1 .\\
% \end{CD}
% $$
% In other words, the group scheme $N$ is the push-out of $\widetilde N$ 
% by the morphism $_2T \to T$ (for the theory of extensions, 
% see~\cite[XVII, appendice I.2]{SGA3}).
% We define then $S \subset \widetilde N$ as the pushout 
% of $\widetilde N$ by the map  $_2T \to {_{n^2}T}$.
% Then $S$ is an extension of $W$ by the finite
% split group of multiplicative type $_{m^2}T$. This completes
% the construction of $S$. By construction, $S$ satisfies conditions
% (1) of Theorem~\ref{thm1}. The proof of (2) is the same as
% in part (a).  \qed
%%%%%%%%%%%%%%%%%%%%%%%%%%%%%%%%%%%

\begin{remark} \label{rem.correction}
In the case where $X = \Spec(K)$ for some field $K/k$,
Theorem~\ref{thm1}(a) reduces to~\cite[Theorem 1.1(a)]{cgr}, and our
proof proceeds along similar lines. Note however, that there is a
small mistake in the proof of~\cite[Theorem 1.1(a)]{cgr}. On page 565
in~\cite{cgr}, in the setting of Lemma~\ref{lem2} above, we said
that the obstruction $\Delta(a)$ (denoted by $\delta([a])$ there) to
lifting $a$ to $H^1(X, S)$ lies in $H^2(X, \mu)$, instead of 
$H^2(X, \, _a\mu)$. This mistake is corrected in the proof of 
Lemma~\ref{lem2} in the present paper. As a consequence, 
the (corrected) argument in
this paper is a bit longer than in~\cite{cgr}, and the group $S$ is
a bit larger; here $S \cap T = {_{n^2}T}$, where as 
in~\cite{cgr} $S \cap T = _nT$.
\end{remark}

\section{Toral torsors and a theorem of Grothendieck}
\label{sect.toral}

Let $X$ be a scheme and $G$ be a smooth affine group scheme over
$X$. Assume that the connected component $G^0$ is reductive. We say
that a $G$-torsor $E$ over $X$ is {\em toral} if the twisted
$X$-group scheme $^EG$ admits a maximal torus defined over $X$. We
denote by $H^1_{toral}(X, G) \subset H^1(X, G)$ the set of toral
classes. The following lemma is well known.

\begin{lem} \label{lem.toral}
Assume that $G^0/X$ admits a maximal
$X$-torus $T$. Then
$$ H^1_{toral}(X, G) =
{\rm Im}\Bigl( H^1(X,N_G(T)) \to
H^1(X,G) \Bigr) \, .  $$
\end{lem}

\begin{proof} Let $E/X$ be a $G$-torsor.
The functor ${\cal T}/X$ of maximal tori of $^EG$ is representable
by a separated smooth scheme $\Sigma$ of finite type over $X$
\cite[XII.1.10]{SGA3}.
In fact, $\Sigma$ is the $E$-twist of homogeneous
space $G/N_G(T)$ (whose points represent maximal tori in $G$);
equivalently, $\Sigma$ can be thought of as the quotient $E/N_G(T)$ (see
\cite[XXIV.4.2.1]{SGA3}). So the following are equivalent:

\begin{enumerate}

\item $^EG$ has a maximal $X$-torus,

\smallskip

\item  ${\cal T}(X) \not = \emptyset$,

\smallskip
         
\item  $(E/N_G(T))(X) \not = \emptyset$.
         
\end{enumerate}
By \cite[III,\,\S 4,\,Prop. 4.6]{DG}, condition (3) is equivalent to
$$       
[E] \in  {\rm Im}\Bigl( H^1(X,N_G(T)) \to H^1(X,G)\Bigr) \, ,
$$       
and the lemma follows.
\end{proof}

The following theorem of Grothendieck tells us that if $G$ 
is a reductive group scheme over a semilocal ring $k$ then 
every $G$-torsor over $k$ is toral. 

\begin{thm} \label{thm3} (\cite[XIV.3.20]{SGA3}).
Let $G$ be a reductive group scheme defined over a semilocal ring
$k$. Then $G$ admits a maximal $k$-torus $T$. \qed
\end{thm}

%%%%%%%%%%%%%%%%%%%%%%%%%
% In the sequel we will use the following
% corollary of this theorem. The special case where $k$ is a field
% is particularly well known and frequently used;
% cf.,~\cite[Lemma III.2.2.1]{serre-gc}.
%%%%%%%%%%%%%%%%%%%%%%%%%
The corollary below will be of particular interest to us in the sequel.

\begin{cor} \label{cor-thm3}
Let $G$ be a smooth affine reductive groups scheme defined over 
a semilocal ring $k$. Suppose $T$ is a maximal $k$-torus of $G$.
Then the natural map $H^1(R,N_G(T)) \to H^1(R,G)$ is
surjective for any semilocal ring $R/k$.
\end{cor}

\begin{proof} By Theorem~\ref{thm3} every $G$-torsor over $\Spec(R)$
is toral. That is, $$H^1(R, G)_{toral} = H^1(R, G) \, . $$
The corollary now follows from Lemma~\ref{lem.toral}.
\end{proof}

\section{Proof of Theorem \ref{thm1}(c)}\label{proof2}

Throughout this section $k$ will denote a semilocal ring and
$G$ an affine connected reductive group scheme 
defined over $k$. Suppose that the radical torus of $G$
is isotrivial. We will now proceed to prove Theorem~\ref{thm1}(c) 
in four steps. 

\medskip

\noindent {\bf Case $1$. $G$ split, semisimple and adjoint.}
That is, $G=G_0\times_{\mathbb{Z}}\,k$, where $G_0$ is an
adjoint split group defined over $\mathbb{Z}$. Let $T_0$ be a maximal 
split torus in $G_0$ defined over $\bbZ$ 
and let $S_0' \subset N_{G_0}(T_0)$ be
the finite subgroup satisfying the conditions of Proposition~\ref{prop0}
constructed in the previous section. Then $S' = S_0' \otimes_{\bbZ} k$
satisfies the same conditions in $G$, relative to 
the maximal torus $T = T_0 \otimes_{\bbZ} k$ of $G$. 
The desired conclusion now follows from Proposition~\ref{prop0}. 
% Taking
% $S=S_0\times_{\mathbb{Z}} \,k$, where $S_0$ is the finite 
% subgroup in $G_0$ constructed in the previous section, we see 
% that in this case Theorem~\ref{thm1}(c) can be proved 
% by exactly the same argument as Theorem~\ref{thm1}(b) 
% in the previous section.

\medskip
\noindent {\bf Case $2$. $G$ is a quasi-split semisimple and adjoint}. 
In this case $G$ is
$k$-isomorphic to the twist ${_{h_*(a)}(G_1)}$, where $G_1$ is a
split adjoint group scheme over $k$ of the same type as $G$,
$$
a\in Z^1_{\acute et}(k,\Aut_{\rm Dyn}\,( D)) \, , $$ and $D$ is the
Dynkin scheme of $G_1$, relative to a maximal split
$k$-torus $T_1 \subset G_1$;
see \S\ref{sect.prel}.  The cocycle
$h_*(a)$ preserves the maximal torus $T_1$ and the finite subgroup 
of $N_{G_1}(T_1)$
constructed in Case $1$. In Case 1 we called this finite subgroup $S'$; 
now we will denote it by $S_1'$. Recall that $S_1'$ satisfies the 
conditions of Proposition~\ref{prop0} relative to 
$T_1$; that is, $S_1'$ normalizes $T_1$ and projects surjectively 
onto $W_1 = N_{G_1}(T_1)/T_1$. Now observe that the group 
$_{h_*(a)}(S_1')$ satisfies the same conditions in $G = 
{_{h_*(a)}(G_1)}$, relative to the maximal $k$-torus 
${_{h_*(a)}(T_1)}$. The desired conclusion now follows from
Proposition~\ref{prop0}.

\medskip

\noindent {\bf Case 3: $G$ is semisimple and adjoint.} 
% Let $G$ be an adjoint
% semisimple $k$-group. 
In this case $G$ is $k$-twisted form of a Chevalley group
$G_0$ \cite[XXIII.5.7]{SGA3}. In other words,
there exists a cocycle $b \in Z^1(k, \Aut(G_0))$ 
such that $G \cong {_bG_0}$. Let $a$ be the image
of $b$ under the projection $$Z^1_{fppf}(k, \Aut(G_0)) \to
Z^1_{fppf}(k, {\rm Aut}\,(D)) \, ,$$ where $D$ is the Dynkin scheme
of $G_0$. Consider the following commutative exact diagram of
pointed sets
$$
\begin{CD}
&& H^1(k,\,  \Aut(G_0)) @>>>  H^1(k, \, \Aut(D)) \\
&& @A{f_{h_*(a)}}AA @A{f_a}AA \\
H^1(k, \, {_{h_*(a)}G_0}) @>>>  H^1(k,\,  {_{h_*(a)}\Aut(G_0)})
@>>>  H^1(k, \, _a\Aut(D)) \\
\end{CD}
$$
where $f_{h_*(a)}$ and $f_a$ stand for  the ``torsion" bijections; 
see \S\ref{sect.prel}.  By a diagram chase, there exists $[c] \in
H^1(k, \, {_{h_*(a)}G_0})$ such that $G$ is isomorphic to the
twisted group $_c\Bigl( {_{h_*(a)}G_0}\Bigr)$, i.e.,~$G$ is a
$k$-inner form of the quasi-split group  ${_{h_*(a)}G_0}$. 

By Case $2$ we know that Theorem~\ref{thm1}(c) holds for
$G_1 = {_{h_*(a)}G_0}$. That is,
there exists a maximal torus $T_1 \subset G_1$ defined 
over $k$ and a finite $k$-subgroup $S_1 \subset N_G(T_1)$, such that
$S_1$ is an extension of a twisted constant group
scheme by a finite $k$-group of multiplicative type and
the natural map $ H^1(k, S_1) \lra H^1(k, N_{G_1}(T_1))$ is surjective.
Moreover, by Corollary~\ref{cor-thm3} the map
$H^1(k, N_{G_1}(T_1)) \to H^1(k, G_1)$ is also surjective.
We conclude that the map $H^1(k,S_1) \to H^1(k,G_1)$ is surjective.
We may thus assume that $c$ takes values in $S_1$. 

Now set $S := \, _cS_1$. Then $S$ embeds in
$_cT_1$. Consider the diagram
% \begin{array}{ccc}
% H^1(X,S) & \xrightarrow{\pi_0} & H^1(X,\,  N_{G_1}(T_1)) \\
%\uparrow &      & \uparrow \\
% H^1(X, {_cS}) & \xrightarrow{\pi} & \; H^1(X,N_G(_cT_1))
% \end{array} 
\[ \begin{CD}
H^1(X,S_1) @>{\pi_0}>> H^1(X,  N_{G_1}(T_1)) \\
@AAA @AAA \\
H^1(X, {_cS_1}) @>{\pi}>> H^1(X, N_G(_cT_1))
\end{CD} \]
where the vertical arrows are the ``torsion"  bijections; 
see~\S\ref{sect.prel}.  Since $\pi_0$ is surjective, so is $\pi$.

\medskip

\noindent {\bf Case $4$. $G$ is reductive and the radical 
torus $C = \operatorname{rad}(G)$ is isotrivial.} 
% Let $G$ be a reductive group scheme over $k$. 
Consider the semisimple $k$-group $H=[G,G]$.
Let $Z$ be the center of $H$,
$G' =G/Z$ and $f:G \to G'$ be the natural projection. Then
$G'\simeq C'\times H'$, where $C'=C/ C \cap H$ and $H'=H/Z$ is an
adjoint semisimple group. Since we are assuming that $C$ is an
isotrivial $k$-torus,
there exists a finite \'etale surjective covering $\widetilde k/ k$
which splits $C$. Note that $\widetilde k$ is a semilocal ring and
$C' \times_k \widetilde k$ is also a split torus.

Let $m$ be the degree of the 
covering $\widetilde k/k$ and let $\mu$ be the $m$-torsion
subgroup of $C'$. Note that the canonical mapping $H^1(k,\mu)\to
H^1(k,C')$ is surjective. Indeed, the restriction-corestriction
formula \cite[0.4]{CTS2}
\[ \times\, m = \operatorname{Cor}_k^{\widetilde k} \circ
\operatorname{Res}_k^{\widetilde k}:\ H^1(k,C')\to H^1(k,C') \]
together with the fact that $H^1(\widetilde k,C')=0$ (Hilbert's
Theorem 90) imply that the map
\[ \times \,m \colon H^1(k,C') \rightarrow H^1(k,C') \]
is trivial.

Let $T'$ and $S'\subset N_{H'}(T')$ be the subgroups constructed in
Case $3$ for $H'$ and let $X/k$ be a scheme satisfying
condition~(\ref{e.C}). Then the canonical morphism
$\pi':H^1(X,\mu\times S') \to H^1(X,N_{G'}(T'))$ is surjective. We
claim that $S=f^{-1}(\mu\times S')$ is as required, i.e.
$H^1(X,S)\to H^1(X,N_G(T))$ is surjective where $T=f^{-1}(C'\times
T')$.

Indeed, the exact sequences $1\to Z\to N_G(T)\to N_{G'}(T') \to 1$
and $1\to Z \to S\to S' \to 1$ give rise to a commutative diagram
$$
% \begin{array}{ccccccc}
% H^1(X,Z) & \xrightarrow{} & H^1(X,N_G(T)) & \xrightarrow{g_1} &
% H^1(X, N_{G'}(T'))
% & \xrightarrow{g_2} & H^2(X,Z) \\
% \uparrow  & &\uparrow \, \pi &    & \uparrow \, \pi'   &  
% & \uparrow \, {id} \\
% H^1(X,Z) & \xrightarrow{} & H^1(X, S) & \xrightarrow{h_1} &
% H^1(X,\mu\times S') \xrightarrow{h_2} & H^2(X,Z)
% \end{array}
\begin{CD}
H^1(X,Z) @>>> H^1(X,N_G(T)) @>{g_1}>> H^1(X, N_{G'}(T')) @>{g_2}>> H^2(X,Z) \\
@AAA  @A{\pi}AA @A{\pi'}AA  @A{\operatorname{id}}AA \\
H^1(X,Z) @>>> H^1(X, S) @>{h_1}>> H^1(X,\mu \times S') @>{h_2}>> H^2(X,Z)
\end{CD}
$$
Here $g_2,h_2$ are connecting homomorphisms \cite[IV.4.3.4]{Gir}.
Fix an element $[a]\in H^1(X,N_G(T))$ and let $[b]=g_1([a])$. Since
$\pi'$ is surjective, there is a class $[c]\in {H^1(X,\mu\times
S')}$ such that $\pi'([c])=[b]$. Since $h_2([c])= g_2 \pi'([c])=0$,
there is $[d]\in H^1(X,S)$ such that $h_1([d])=[c]$. Thus the
classes $[a]$ and $\pi ([d])$ have the same image in
$H^1(X,N_{G'}(T'))$. A twisting argument shows that the map $H^1(X,
{_dZ})\to g^{-1}_1(g_1([a]))$ is surjective;
see~\cite[III.3.2.4]{Gir}. Since $Z\subset S$, we have $_dZ\subset
{_dS}$ implying $[a]\in {\rm Im}\, \pi$. This completes the proof of
Theorem~\ref{thm1}(c). \qed

\section{Proof of Theorem~\ref{thm2} and an application}
\label{sect.thm2}

Theorem~\ref{thm2} stated in the introduction is an easy consequence
of Theorem~\ref{thm1}(c) and Corollary~\ref{cor-thm3}. Indeed,
choose $T$ and $S$ as in Theorem~\ref{thm1}(c) and let $R$ be a semilocal
ring containing $k$. Since $X = \Spec(R)$ satisfies
condition~\eqref{e.C}, Theorem~\ref{thm1}(c) tells us that the
natural map $H^1(R, S) \to H^1(R, N(T))$ is surjective. By
Corollary~\ref{cor-thm3} the map $H^1(R,N_G(T)) \to H^1(R,G)$ is
also surjective, and Theorem~\ref{thm2} follows.
\qed

We will now discuss an application of Theorem~\ref{thm2}.
Let $G$ be a linear algebraic group defined over a field $k$
and $\pi \colon Y \to X$ be a $G$-torsor over a $k$-scheme $X$.
As usual, we will say that $\pi$ admits
reduction of structure to a $k$-subgroup $S \subset G$
if the class in $H^1(X, G)$ represented by $\pi$ lies in the image
of the natural map $H^1(X, S) \to H^1(X, G)$. Equivalently,
$\pi$ admits reduction of structure to $S$ if
there exists a $G$-equivariant morphism $Y \to G_X/S_X$.

Suppose $U \to X$ is a morphism and $Y_U$
is the pull-back of $Y$ to $U$:
$$
\begin{CD}
Y_U  @>>>  Y \\
@V{\pi_U}VV @V{\pi}VV \\
U @>>> X \, . \\
\end{CD}
$$
We say that $Y$ admits reduction of structure to $S$ {\em over
$U$} if the $G$-torsor $\pi_{U} \colon Y_U \to U$
admits reduction of structure to $S$. 

\begin{prop} \label{prop.section}
Let $G$ be a linear algebraic group defined over a field $k$.
Assume that the connected component $G^0$ is reductive and
either $G$ is connected or $k$ is algebraically closed.
Then there exists a $k$-subgroup $S \subset G$ with the following
property.  For any $G$-torsor $\pi \colon Y \to X$ over
an affine $k$-scheme $X$
and any finite collection of (not necessarily closed) points
$x_1, \dots, x_n \in X$,
there exists an open subscheme $U \subset X$ containing
$x_1, \dots, x_n$ such that $\pi$ admits reduction
of structure to $S$ over $U$.
\end{prop}

\begin{proof} Let $R = \cal{O}_{x_1, \dots, x_n}$ be the semilocal
ring of $X$ at $x_1, \dots, x_n$.  By Theorem~\ref{thm2} $\pi$
admits reduction of structure to $S$ over $\Spec(R) \subset X$. That
is, there exists a $G$-equivariant morphism $\phi \colon Y_R \to
G_R/S_R$.

Since $R$ is, by definition, the direct limit of
$\cal{O}_X(U)$, as $U$ ranges over the open subsets of $X$
containing $x_1, \dots, x_n$, $\phi$ extends over
some open subscheme $X_0$ of $X$ containing $x_1, \dots, x_n$.
In other words, $\pi$ admits reduction of structure to $S$ over $X_0$.
\end{proof}

\section{Torsors on affine spaces}
\label{sect.affine}

Let $k$ be a field of characteristic $0$,
and $\overline{k}$ be an algebraic closure of $k$.
In this section we will apply Theorem~\ref{thm1}(c)
in the case where $X$ is the affine space $\bbA_k^{\; n}$.  
The key observation here is that
$\bbA_{\overline{k}}^{\; n}$ is simply connected.
By the fundamental exact sequence for $\pi_1$~\cite[IX.6.1]{SGA1},
we have an isomorphism 
$$\pi_1( \bbA_k^{\; n},\overline 0) \simlgr \Gal(\overline{k}/k),$$ where
$\pi_1( \bbA_k^{\; n},\overline 0)$ stands for the algebraic 
fundamental group of $ \bbA_k^{\; n}$ relative to the base point 
$\overline 0 : \Spec(\overline k) \to  \bbA_k^{\; n}$.
In other words, every finite
\'etale cover of $\bbA_{k}^{\; n}$ is of the form
$\bbA_{K}^{\; n}$, where $K/k$ is an \'etale $k$-algebra.
Since $\Pic(\bbA_{K}^{\; n} )=0$, this implies that 
$X=\bbA_k^{\; n}$ satisfies condition~\eqref{e.C}.

\begin{prop}\label{prop.affine} Let $k$ be a field
of characteristic zero, $n \geq 0$ be an integer, and
$G$ be a (connected) reductive group over $k$.  Then
$$
H^1(k,G) \simlgr H^1(\bbA_{k}^{\; n},G)_{toral}.
$$
In other words, every toral torsor on $\bbA_{k}^{\; n}$
is constant.
\end{prop}

\begin{proof}
Since $k$ is a field, Theorem~\ref{thm1}(c)
applies to $G$.  Let $S \subset G$ be the finite $k$-subgroup
as in Theorem \ref{thm1}.
As we noted before the statement of the proposition, 
$X=\bbA_k^{\; n}$ satisfies condition~\eqref{e.C}.
Thus the natural map
$ H^1( \bbA_{k}^{\; n}, S) \to H^1( \bbA_{k}^{\; n}, N(T))$ 
is surjective.  By Lemma~\ref{lem.toral}, the map
$$
H^1( \bbA_{k}^{\; n}, S) \to H^1( \bbA_{k}^{\; n},G)_{toral}
$$
is also surjective. The kernel of the natural map
$H^1(\bbA_k^{\; n},S) \to H^1(\bbA_{\overline{k}}^{\; n},S)$ consists
of those $S$-torsors on $\bbA_k^{\; n}$ which become trivial
on $\bbA_{\overline{k}}^{\; n}$.
Since $\bbA_{\overline{k}} \to \bbA_k$ is a Galois cover, with
Galois group $\Gal(\overline{k}/k)$, this kernel is
$H^1(k, S(\bbA_{\overline{k}}^{\; n}))$, where $H^1$ 
stands for Galois cohomology. Since
$S( \bbA_{\overline{k}}^{\; n}) = S(\overline{k})$,
this yields an exact sequence
$$
1 \to H^1(k, S(\overline{k}))
\to  H^1(\bbA_k^{\; n},S ) \to H^1(\bbA_{\overline{k}}^{\; n},S ) .
$$
Since $\bbA_{\overline{k}}^{\; n}$ is simply connected,
$H^1(\bbA_{\overline{k}}^{\; n},S )=1$, and hence
the map 
\[ H^1(k, S(\overline{k})) \to  H^1(\bbA_k^{\; n},S) \]
is surjective.
The commutative exact diagram of pointed sets
$$
\begin{CD}
 H^1(k, S(\overline{k})) @>>>   H^1(\bbA_k^{\; n},S ) @>>>1 \\
@VVV @VVV \\
 H^1(k, G) @>>>   H^1(\bbA_k^{\; n},G )_{toral} \\
@VVV @VVV \\
1 && 1 \\
\end{CD}
$$
shows that the natural map $H^1(k, G) \to H^1(\bbA_k^{\; n},G)_{toral}$
is surjective.  This map is also injective. Indeed, suppose
$G$-torsors $T_1 \to \Spec(k)$ and $T_2 \to \Spec(k)$ map to the
same $G$-torsor $Y \to \bbA_k^{\; n}$, i.e., $Y \simeq T_i
\times_{\Spec(k)} \bbA_k^{\; n}$ for $i = 1, 2$. Then both $T_1$ and
$T_2$ are isomorphic to the fiber of $Y$ over $0 \in \bbA_k^{\; n}$.
Hence, $T_1$ and $T_2$ represent the same class in $H^1(k, G)$. We
conclude that the map $H^1(k, G) \to H^1(\bbA_k^{\; n},G)_{toral}$ is an
isomorphism.
\end{proof}

\begin{remark} \label{rem.OS}
There are examples of non constant $G$-torsors $P$ over
affine spaces; see Ojanguren-Sridharan~\cite{OS}
(cf.~also \cite[VII.10]{K}).  Proposition~\ref{prop.affine} tells us
that in these examples the twisted groups $^P G$
do not carry maximal tori.
\end{remark}

\begin{remark} \label{rem.loop}
As we pointed out in the introduction, the scheme 
\[ X = \Spec(k[x_1^{\pm 1}, \dots, x_n^{\pm 1}]) \]
also satisfies
condition~\eqref{e.C} (in characteristic zero), so in this case the
map $H^1(X, S) \to H^1(X, G)_{toral}$ is also surjective. This fact
is used in~\cite{GP2}.
\end{remark}

\end{document}